\documentclass[a4j]{article}

\font\es=eufm10

\def\R{\mbox{\boldmath $R$}}
\def\C{\mbox{\boldmath $C$}}
\def\gP{\mbox{\es {P}}}

\def\J{\mbox{\es {J}}}
\def\CC{\mbox{\es {C}}}
\def\Z{\mbox{\boldmath $Z$}}
\def\H{\mbox{\boldmath $H$}}
\def\be{\mbox{\boldmath $e$}}
\def\b1{\mbox{\boldmath $1$}}
\def\bi{\mbox{\boldmath $i$}}
\def\bj{\mbox{\boldmath $j$}}
\def\bk{\mbox{\boldmath $k$}}
\begin{document}
\title { {\bf Spinor-generators of compact exceptional  \\
        Lie groups $F_4, E_6$ and $E_7$ }}

\author {Takashi Miyasaka, Osamu Shukuzawa and Ichiro Yokota}
\date{}
\maketitle
\section {Introduction} 

\hspace*{5pt} We know that any element $A$ of the group $SO(3)$ can be represented as
$$
         A = A_1A_2{A_1}', \quad A_1, {A_1}' \in SO_1(2), A_2 \in SO_2(2) $$
\noindent where $SO_k(2) = \{A \in SO(3) \, | \, A\be_k = \be_k \} (k = 1, 2)$ ([1]). In the present paper, we shall show firstly that the similar results hold for the groups $SU(3),$ and $Sp(3)$ (Theorem 1). Secondly, we shall  show that any element $\alpha$ of the simply connected compact Lie group $F_4$ (resp. $E_6$ ) can be represented as 
$$
          \alpha = \alpha_1\alpha_2{\alpha_1}',\quad
        \begin{array} {c}
     \alpha_1, {\alpha_1}' \in Spin_1(9),  \alpha_2 \in Spin_2(9) \cr
  (\mbox{resp.} \alpha_1, {\alpha_1}' \in Spin_1(10), \alpha_2 \in Spin_2(10)) \end{array} $$     
\noindent where $Spin_k(9) = \{\alpha \in F_4 \, | \, \alpha E_k = E_k \}$ (resp. $Spin_k(10) = \{\alpha \in E_6 \, | \, \alpha E_k = E_k \}$ (Theorem 5 (resp. Theorem 7))). Lastly, we shall show that any element $\alpha$ of the simply connected compact Lie group $E_7$ can be represented as
$$
   \alpha = \alpha_1\alpha_2{\alpha_1}'{\alpha_2}'{\alpha_1}'', \quad
       \alpha_1, {\alpha_1}', {\alpha_1}'' \in Spin_1(12),
       \alpha_2, {\alpha_2}'  \in Spin_2(12) $$
\noindent where $Spin_k(12) = \{ \alpha \in E_7 \, | \, \alpha\kappa_k = \kappa_k\alpha, \alpha\mu_k = \mu_k\alpha \}$ (Theorem 10).
\vspace{1mm}

\hspace*{5pt} In this paper we follow the notation of [2].
\section {Spinor-generators of the groups $SO(3), SU(3)$ and $Sp(3)$ }
\hspace*{5pt} Let $\H$ be the quaternion field with basis $1$, $\bi, \bj$ and $\bk$ over $\R$. Then we can express each element $a = a_0 + a_1\bi + a_2\bj + a_3\bk \in \H$ in the following polar form
$$
    a = r(\cos\theta + u\sin\theta),\quad
    u^2 = -1 \, (u \in \H), \, \, r = |a| = \sqrt{\sum_{k=0}^{3}{a_k}^2}, \theta \in \R. $$
\noindent Hereafter, we briefly denote by $re^{u\theta}$ an element $r(\cos\theta + u\sin\theta)$ after the model of complex numbers.
\vspace{1mm}

\hspace*{5pt} The classical groups $SO(n), SU(n)$ and $Sp(n)$ are respectively defined by
$$ 
    \begin{array}{l}
       SO(n) = \{ A \in M(n, \R) \, | \, {}^tAA = E, \mbox{det}A = 1 \}, \\
       SU(n) = \{ A \in M(n, \C) \, | \, A^*A = E, \mbox{det}A = 1 \}, \\
       Sp(n) = \{ A \in M(n, \H) \, | \, A^*A = E \}
    \end{array} $$
\noindent where we follow the usual convention for matrices : $M(n, K)$ ( = the set of square matrices of order $n$ with coefficients in $K = \R, {\bf C}$ or $\H$ ), ${}^tA, A^* ( = \overline{{}^tA}), E$ ( = the unit matrix) and det ( = the determinant).
\vspace{2mm}

\hspace*{5pt} {\bf Theorem 1.} (1) {\it Any element $A \in SO(3)$ can be represented as}
$$ 
      A = A_1A_2{A_1}', \quad A_1, {A_1}' \in SO_1(2), A_2 \in SO_2(2) $$
\noindent {\it where} $SO_k(2) = \{ A \in SO(3) \, | \, A\be_k = \be_k \} \cong Spin(2) \, ( k = 1, 2), \be_1 = {}^t(1, 0, 0), \be_2 = {}^t(0, 1, 0)$.
\vspace{1mm}

\hspace*{5pt} (2) {\it Any element $A \in SU(3)$ can be represented as}
$$ 
      A = A_1A_2{A_1}', \quad A_1, {A_1}' \in SU_1(2), A_2 \in SU_2(2) $$
\noindent {\it where} $SU_k(2) = \{ A \in SU(3) \, | \, A\be_k = \be_k \} \cong Spin(3) \, ( k = 1, 2)$.
\vspace{1mm}

\hspace*{5pt} (3) {\it Any element $A \in Sp(3)$ can be represented as}
$$ 
      A = A_1A_2{A_1}', \quad A_1, {A_1}' \in Sp_1(2), A_2 \in Sp_2(2) $$
\noindent {\it where} $Sp_k(2) = \{ A \in Sp(3) \, | \, A\be_k = \be_k \} \cong Spin(3) \, ( k = 1, 2)$.
\vspace{2mm}

\hspace*{5pt} {\bf Proof} It suffices to prove (3), because we can reduce (1) and (2) to the particular case of (3) in the proof below. First, for a given element $A \in Sp(3)$, suppose $A\be_1 = {}^t(a_1, a_2, a_3), a_2 \neq 0 \, (a_k \in \H \, ( k = 1, 2, 3))$. Then there exists an element $u \in \H$ satisfying $u^2 = -1$ and a real number $\alpha \in \R$ such that $a_3a_2^{-1} = (|a_3|/|a_2|)e^{u\alpha}$. Choose $\theta \in \R$ such that $\cot\theta = |a_3|/|a_2|$ and set
$$ 
        B_1 = \pmatrix{1 & 0 & 0 \cr
               0 & e^{u\alpha/2}\cos \theta & - e^{-u\alpha/2}\sin \theta \cr
               0 & e^{u\alpha/2}\sin \theta &  e^{-u\alpha/2}\cos \theta }
             \in Sp_1(2). $$
\noindent Then we get 
$$
         B_1A\be_1 = {}^t(b_1, 0, b_3), \quad
         b_1, b_3 \in \H. $$
\noindent Next suppose $b_3 \neq 0.$ Then there exists an element $v \in \H$ satisfying $v^2 = -1$ and a real number $\beta \in \R$ such that $b_1{b_3}^{-1} = (|b_1|/|b_3|)e^{v\beta}$. Choose $\varphi \in \R$ such that $\cot\varphi = - |b_1|/|b_3|$ and set
$$
    B_2 = \pmatrix{e^{-v\beta/2}\cos\varphi & 0 & - e^{v\beta/2}\sin\varphi \cr
                     0 & 1 & 0 \cr
                    e^{-v\beta/2}\sin\varphi & 0 & e^{v\beta/2}\cos\varphi}
                   \in Sp_2(2). $$ 
\noindent Then we get
$$
            B_2B_1A\be_1 = {}^t(c_1, 0, 0), \quad c_1 \in \H. $$
\noindent Since $|c_1| = 1,$ we can say $c_1 = e^{w\gamma} \, (w^2 = -1, w \in \H, \gamma \in \R).$ Set
$$
            {B_2}' = \pmatrix{e^{-w\gamma} & 0 & 0 \cr
                              0 & 1 & 0 \cr
                              0 & 0 & e^{w\gamma} } \in Sp_2(2). $$
\noindent Then, since it follows ${B_2}'B_2B_1A\be_1 = \be_1,$ {\it i.e.}, ${B_2}'B_2B_1A \in Sp_1(2)$, we can set ${B_2}'B_2B_1A = {B_1}' \in Sp_1(2)$. This implies
$$
     A = A_1A_2{A_1}', \quad A_1, {A_1}' \in Sp_1(2), \, A_2 \in Sp_2(2). $$
\section { Some elements of $Spin_k(9), Spin_k(10)$ and $Spin_k(12).$ }   
\hspace*{5pt} As for the definitions of $Spin_k(9), Spin_k(10)$ and $Spin_k(12) (k = 1, 2)$, see Sections 4, 5 and 6.
\vspace{3mm}

\hspace*{5pt} {\bf Lemma 2} (Section 4 and [2]). (1) {\it Let $\alpha_1(a)$ be the mapping $\alpha(a)$ defined in} [2] {\it Lemma} 2.(1). {\it Then $\alpha_1(a)$ belongs to $Spin_1(9) \subset Spin_1(10) \subset Spin_1(12)$.}
\vspace{1mm}

\hspace*{5pt} (2) {\it For $a \in \CC, a \neq 0,$ let $\alpha_2(a) : \J \to \J$ be the mapping defined by changing all of the indices from $k$ to $k + 1$} ({\it index modulo} 3) {\it in the definition of $\alpha(a)$ of} [2] {\it Lemma} 2.(1). {\it Then $\alpha_2(a)$ belongs to $Spin_2(9) \subset Spin_2(10) \subset Spin_2(12)$.}
\vspace{3mm}

\hspace*{5pt} {\bf Lemma 3} (Section 5 and [2]). (1) {\it Let $\beta_1(a)$ be the mapping $\beta(a)$ defined in} [2] {\it Lemma} 2.(2). {\it Then $\beta_1(a)$ belongs to $Spin_1(10) \subset Spin_1(12)$.}
\vspace{1mm}

\hspace*{5pt} (2) {\it For $a \in \CC, a \neq 0,$ let $\beta_2(a) : \J^C \to \J^C$ be the mapping defined by changing all of the indices from $k$ to $k + 1$} ({\it index modulo} 3) {\it in the definition of $\beta(a)$ of} [2] {\it Lemma} 2.(2). {\it Then $\beta_2(a)$ belongs to $Spin_2(10) \subset Spin_2(12)$.}
\vspace{3mm}
  
\hspace*{5pt} {\bf Lemma 4} (Section 6 and [2]). (1) {\it Let $\gamma_1(a)$ be the mapping $\gamma(a)$ defined in} [2] {\it Lemma} 3.(1). {\it Then $\gamma_1(a)$ belongs to  $Spin_1(12)$.}
\vspace{1mm}

\hspace*{5pt} (2) {\it For $a \in \CC, a \neq 0,$ let $\gamma_2(a) : \gP^C \to \gP^C$ be the mapping defined by changing all of the indices from $k$ to $k + 1$} ({\it index modulo} 3) {\it in the definition of $\gamma(a)$ of} [2] {\it Lemma} 3.(1). {\it Then $\gamma_2(a)$ belongs to $Spin_2(12)$.}

\hspace*{5pt} (3) {\it Let $\delta_1(a)$ be the mapping $\delta(a)$ defined in} [2] {\it Lemma} 3.(2). {\it Then $\delta_1(a)$ belongs to  $Spin_1(12)$.}
\vspace{1mm}

\hspace*{5pt} (4) {\it For $a \in \CC, a \neq 0,$ let $\delta_2(a) : \gP^C \to \gP^C$ be the mapping defined by changing all of the indices from $k$ to $k + 1$} ({\it index modulo} 3) {\it in the definition of $\delta(a)$ of} [2] {\it Lemma} 3.(2). {\it Then $\delta_2(a)$ belongs to $Spin_2(12)$.}
\section {$Spin(9)$-generatos of the group $F_4$}
\hspace*{5pt} The simply connected compact Lie group $F_4$ is given by
$$
    F_4 = \{ \alpha \in \mbox{Iso}_R(\J) \, | \, \alpha(X \times Y)
        = \alpha X \times \alpha Y \}. $$
\noindent The group $F_4$ has subgroups
$$
    Spin_k(9) = \{ \alpha \in F_4 \, | \, \alpha E_k = E_k \} \quad (k = 1, 2), $$
\noindent where $E_1 = (1, 0, 0; 0, 0, 0), E_2 = (0, 1, 0; 0, 0, 0) \in \J$, which is isomorphic to the usual spinor group $Spin(9)$ ([2],[3]).
\vspace{3mm}

\hspace*{5pt} {\bf Theorem 5.} {\it Any element $\alpha \in F_4$ can be represented as}
$$
    \alpha = \alpha_1\alpha_2{\alpha_1}', \quad 
       \alpha_1, {\alpha_1}' \in Spin_1(9),\alpha_2 \in Spin_2(9). $$

\hspace*{5pt} {\bf Proof} For a given element $\alpha \in F_4$, it suffices to show that there exist $\alpha_1 \in Spin_1(9)$ and $\alpha_2 \in Spin_2(9)$ such that $\alpha_2\alpha_1\alpha E_1 = E_1$. Now, for $\alpha E_1 = (\xi_1, \xi_2, \xi_3; x_1, x_2, x_3) = X_0$, choose $a \in \CC$ such that $(a, x_1) = 0, |a| = \pi/4$, and define $\alpha_1(a) \in Spin_1(9)$ of Lemma 2.(1). Then we get
$$
     \alpha_1(a)X_0 = ({\xi_1}', {\xi_2}', {\xi_3}'; {x_1}', {x_2}', {x_3}') 
                    = X_1. \quad
    {\xi_1}' = \xi_1, {\xi_2}' = {\xi_3}' \in \R, {x_3}' = x_3 \in \CC. $$ 
\noindent If ${x_1}' \neq 0$, define $\alpha_1(\pi{x_1}'/4|{x_1}'|) \in Spin_1(9)$. Then we get
$$
               \alpha_1(\pi{x_1}'/4|{x_1}'|)X_1
          = ({\xi_1}'', {\xi_2}'', {\xi_3}''; 0, {x_2}'', {x_3}'') 
                    = X_2, \quad
        {\xi_1}'' = {\xi_1}', {\xi_k}'' \in \R, {x_k}'' \in \CC. $$ 
\noindent From the condition $X_2 \times X_2 = 0$, we have
$$
   \begin{array}{lll}
 {\xi_2}''{\xi_3}'' = 0, & \quad {\xi_3}''{\xi_1}'' ={x_2}''\overline{{x_2}''},         & \quad {\xi_1}''{\xi_2}'' ={x_3}''\overline{{x_3}''}, \\
    \overline{{x_2}''}\,\overline{{x_3}''} = 0, & \quad {\xi_2}''{x_2}'' = 0,              & \quad {\xi_3}''{x_3}'' = 0.
   \end{array} $$
\noindent Here we distinguish the following cases :
\vspace{1mm}

\hspace*{5pt} (I) When ${\xi_2}'' = 0$. From $ X_2 \times X_2 = 0$, we have ${x_3}''\overline{{x_3}''} = 0$, hence ${x_3}'' = 0$.  Therefore $X_2$ is of the form
$$
        X_2 = ({\xi_1}'', 0, {\xi_3}''; 0, {x_2}'', 0), \quad
        {\xi_1}'' = {\xi_1}', {\xi_3}'' \in \R, {x_2}'' \in \CC. $$
\noindent Choose $b \in \CC$ such that $(b, {x_2}'') = 0, |b| = \pi/4$, and define $\alpha_2(b) \in Spin_2(9)$ of Lemma 2.(2). Then 
$$
       \alpha_2(b)X_2 = (\xi_1^{(3)}, 0, \xi_3^{(3)}; 0, x_2^{(3)}, 0) 
                      = X_3, \quad
       \xi_1^{(3)} = \xi_3^{(3)} \in \R, x_2^{(3)} \in \CC. $$
\noindent If $x_2^{(3)} \neq 0,$ define $\alpha_2(\pi x_2^{(3)}/4|x_2^{(3)}|) \in Spin_2(9)$. Then 
$$
      \alpha_2(\pi x_2^{(3)}/4|x_2^{(3)}|)X_3  
            = (\xi_1^{(4)}, 0, \xi_3^{(4)}; 0, 0, 0) 
                      = X_4, \quad
         \xi_1^{(4)}, \xi_3^{(4)} \in \R. $$
\noindent From $X_4 \times X_4 = 0$, we have $\xi_1^{(4)}\xi_3^{(4)} = 0$. If $\xi_1^{(4)} = 0$, define $\alpha_2(\pi/2) \in Spin_2(9)$. Then
$$
         \alpha_2(\pi/2)X_4 = (\xi_1^{(5)}, 0, 0; 0, 0, 0) 
                            = X_5, \quad
         \xi_1^{(5)} = \xi_3^{(4)} \in \R. $$
\noindent If $\xi_3^{(4)} = 0, X_4$ is also of the form $X_5$. Thus we obtain $X_5 = E_1$, since $\xi_1^{(5)} = \mbox{tr}(X_5) = \mbox{tr}(E_1) = 1.$
\vspace{1mm}

\hspace*{5pt} (II) When ${\xi_2}'' \neq 0$. From $ X_2 \times X_2 = 0$, we have ${\xi_2}''{x_2}'' = 0$, hence ${x_2}'' = 0$. Therefore $X_2$ is of the form 
$$
        X_2 = ({\xi_1}'', {\xi_2}'', 0; 0, 0, {x_3}''), \quad
           {\xi_k}'' \in \R, {x_3}'' \in \CC.  $$
\noindent By considering $\alpha_1(\pi/2)X_2,$ this can be reduced to the case (I).

\hspace*{5pt} We have just completed the proof of Theorem 5.
\section {$Spin(10)$-generators of the group $E_6$}
\hspace*{5pt} The simply connected compact Lie group $E_6$ is given by
$$
     E_6 = \{ \alpha \in \mbox{Iso}_C(\J^C) \, | \, 
          \alpha X \times \alpha Y = \tau\alpha\tau(X \times Y),
          <\alpha X, \alpha Y> = <X, Y> \}. $$
\noindent The group $E_6$ has subgroups 
$$
      Spin_k(10) = \{ \alpha \in E_6 \, | \, \alpha E_k = E_k \} 
                \quad (k = 1, 2), $$
\noindent which is isomorphic to the usual spinor group $Spin(10)$ ([2],[3]).
\vspace{3mm}

\hspace*{5pt} {\bf Lemma 6.} (1) {\it For any element}
$$
          X = (\xi_1, \xi_2, \xi_3; x_1, 0, 0), 
               \quad \xi_k \in C, x_1 \in \CC^C$$
\noindent {\it of $\J^C$, there exists some element} $\alpha_1 \in Spin_1(10)$ {\it such that}$$
          \alpha_1X = ({\xi_1}', {\xi_2}', {\xi_3}'; 0, 0, 0), 
               \quad {\xi_1}' = \xi_1, {\xi_k}' \in C. $$
\hspace*{5pt} (2)  {\it For any element}
$$
          X = (\xi_1, 0, 0; 0, x_2, x_3), 
               \quad \xi_1 \in C, x_k \in \CC^C$$
\noindent {\it of $\J^C$, there exists some element} $\alpha_1 \in Spin_1(9)$ {\it such that}$$
          \alpha_1X = (\xi_1, 0, 0; 0, {x_2}', {x_3}'), 
          \quad {\xi_1}' = \xi_1 \in C, {x_2}' \in \CC^C, {x_3}' \in \CC. $$ 

\hspace*{5pt} {\bf Proof} (1) For $x_1 = p + iq \, (p, q \in \CC)$, choose $a \in \CC , a \neq 0,$ such that $(a, p) = (a, q) = 0$, and define $\alpha_1(\pi a/4|a|) \in Spin_1(9)$ of Lemma 2.(1). Then 
$$
             \alpha_1(\pi a/4|a|)X
               = ({\xi_1}', {\xi_2}', {\xi_2}'={\xi_3}', {\xi_3}'; {x_1}', 0, 0) 
               = X_1.   \quad 
       {\xi_1}' = \xi_1, {\xi_2}' = {\xi_3}' \in C, {x_1}' \in \CC^C.$$ 
\noindent Next, for ${x_1}' = p' + iq' \, (p', q' \in \CC)$, choose $b \in \CC, b \neq 0$, such that $(b, p') = (b, q') = 0$, and define $\beta_1(\pi b/4|b|) \in Spin_1(10)$ of Lemma 3.(1). Then 
$$
             \beta_1(\pi b/4|b|)X_1
                    = ({\xi_1}'', 0, 0; {x_1}'', 0, 0) 
                    = X_2, \quad 
          {\xi_1}'' = \xi_1 \in C, {x_1}'' \in \CC^C. $$  
\noindent Next, for ${x_1}'' = p'' + iq'' \, (p'', q'' \in \CC)$, if $q''\neq 0$, define $\alpha_1(\pi q''/4|q''|) \in Spin_1(9)$ . Then 
$$
             \alpha_1(\pi q''/4|q''|)X_2
                = (\xi_1^{(3)}, \xi_2^{(3)}, \xi_3^{(3)}; p^{(3)}, 0, 0) 
                = X_3, \quad 
        \xi_1^{(3)} = \xi_1, \xi_3^{(3)}=-\xi_2^{(3)} \in C, p^{(3)} \in \CC. $$
\noindent Finally, if $p^{(3)} \neq 0$, define $\beta_1(\pi p^{(3)}/4|p^{(3)}|) \in Spin_1(10)$. Then we get  
$$
             \beta_1(\pi p^{(3)}/4|p^{(3)}|)X_3
                = (\xi_1^{(4)}, \xi_2^{(4)}, \xi_3^{(4)}; 0, 0, 0),                \quad \xi_1^{(4)} = \xi_1, \xi_k^{(4)} \in C $$
\noindent as desired.
\vspace{1mm}

\hspace*{5pt} (2) At first, we show that for any element 
  $$
  Z=(\zeta_1, 0 , 0 ; 0 , z_2 , z_3 ), \quad \zeta_1 \in \R, z_k \in \CC,$$
  \noindent there exists $\alpha_1 \in Spin_1(9)$ such that
  $$
  \alpha_1 Z=(\zeta_1', 0 , 0 ; 0 , z_2' , 0 ), \quad \zeta_1' \in \R, z_2'\in \CC.$$
  \noindent In fact,if $z_2z_3 \neq 0$, choose $t > 0$ such that $\cot(t|z_2z_3|) = - |z_2|/|z_3|$, and define $\alpha_1(t\overline{z_2z_3}) \in Spin_1(9)$. Then we get ($z_3$-part of $\alpha_1(t\overline{z_2z_3})Z) = 0$. If $z_2 = 0$, then $\alpha_1(\pi/2)Z$ is of the form as desired.Now for a given element $X=(\xi_1, 0 , 0 ; 0 , x_2 , x_3 ) \in \J^C$ , express $X=Y+iZ, Y,Z \in \J$ and apply the result above to $Z$, then we get the required form $\alpha_1X=\alpha_1Y+i\alpha_1Z$.
\vspace{3mm}

\hspace*{5pt} {\bf Theorem 7.} {\it Any element $\alpha \in E_6$ can be represented as}
$$
    \alpha = \alpha_1\alpha_2{\alpha_1}', \quad
       \alpha_1, {\alpha_1}' \in Spin_1(10),\alpha_2 \in Spin_2(10). $$

\hspace*{5pt} {\bf Proof} For a given element $\alpha \in E_6$, set $\alpha E_1 = (\xi_1, \xi_2, \xi_3; x_1, x_2, x_3) = X_0 \in \J^C.$  By Lemma 6.(1), we can take $\alpha_1 \in Spin_1(10)$ such that
$$
         \alpha_1X_0
               = ({\xi_1}', {\xi_2}', {\xi_3}'; 0, {x_2}', {x_3}') 
               = X_1, \quad
          {\xi_1}' = \xi_1,$$
\noindent noting that the subspaces $\{(\xi_1,\xi_2,\xi_3;x_1,0,0) \in \J^C \}$ and $\{(0,0,0;0,x_2,x_3) \in \J^C \}$ are invariant under the action of the elements of $Spin_1(10)$, respectively.
From the condition $X_1 \times X_1 = 0$, we have ${\xi_2}'{\xi_3}' = 0$. As a result, the argumant is divided into the following three cases :
\vspace{1mm}

\hspace*{5pt} (I) Case ${\xi_2}' = 0, {\xi_3}'\neq 0$. From $X_1 \times X_1 = 0$, we have ${\xi_3}'{x_3}' = 0$, hence ${x_3}' = 0$. Therefore $X_1$ is of the form
$$
        X_1 = ({\xi_1}', 0, {\xi_3}'; 0, {x_2}', 0), \quad
             {\xi_1}' = \xi_1, $$
\noindent Thus, for $X_1 \in \J^C$, we can take $\alpha_2 \in Spin_2(10)$ such that 
$$
     \alpha_2X_1 = ({\xi_1}'', 0, {\xi_3}''; 0, 0, 0)
                 = X_2, $$
\noindent in the same way as in Lemma 6.(1). Then, from $X_2 \times X_2 = 0$, we have ${\xi_1}''{\xi_3}'' = 0$, tha is,
$$
     X_2 = ({\xi_1}'', 0, 0; 0, 0, 0), \, \, (\tau{\xi_1}''){\xi_1}'' = 1 
         \quad \mbox{or} \quad 
     X_2 = (0, 0,{\xi_3}''; 0, 0, 0), \, \, (\tau{\xi_3}''){\xi_3}'' = 1. $$
\noindent Thus we obtain that there exist some elements $\varepsilon_2(t) \in Spin_2(10)$ and $\alpha_2(\pi/2) \in Spin_2(9)$ such that
$$
    \varepsilon_2(t)X_2 = E_1 \quad \mbox{or} \quad \varepsilon_2(t)\alpha_2(\pi/2)X_2 = E_1, $$
\noindent where $\varepsilon_2(t) \in Spin_2(10)$ is defined by
$$
       \varepsilon_2(t)(\xi_1, \xi_2, \xi_3; x_1, x_2, x_3) = 
     (e^{it}\xi_1, \xi_2, e^{-it}\xi_3; e^{-it/2}x_1, x_2, e^{it/2}x_3),
        \quad t \in \R $$
\noindent (cf. [2] Lemma 10.(1)).
\vspace{1mm}

\hspace*{5pt} (II) Case ${\xi_2}' \neq 0, {\xi_3}' = 0$. From $X_1 \times X_1 = 0$, we have ${\xi_2}'{x_2}' = 0$, hence ${x_2}' = 0$. Therefore $X_1$ is of the form
$$
        X_1 = ({\xi_1}', {\xi_2}', 0; 0, 0,{x_3}'), \quad
         {\xi_1}' = \xi_1.  $$
\noindent Thus, by considering $\alpha_1(\pi/2)X_1$, where $\alpha_1(\pi/2) \in Spin_1(9)$, this can be reduced to Case (I).
\vspace{1mm}

\hspace*{5pt} (III) Case ${\xi_2}' = {\xi_3}' = 0$. By Lemma 6.(2),we can take ${\alpha_1}' \in Spin_1(9)$ such that
$$
      {\alpha_1}'X_1 = ({\xi_1}'', 0, 0; 0, {x_2}'', {x_3}'')
                 = X_2,  \quad 
          {\xi_1}'' = \xi_1, {x_2}'' \in \CC^C, {x_3}'' \in \CC. $$
\noindent Then, from $X_2 \times X_2 = 0$ we have ${x_3}''\overline{x_3''} = 0,$  hence  ${x_3}'' = 0$. Thus, for $X_2 = ({\xi_1}'', 0, 0; 0, {x_2}'', 0) \in \J^C$, we can take $\alpha_2 \in Spin_2(10)$ such that 
$$
      \alpha_2X_2 = (\xi_1^{(3)}, 0, \xi_3^{(3)}; 0, 0, 0) = X_3, $$
\noindent in the same way as in Lemma 6.(1). Hence in a similar way as in Case (I), we obtain the required result.

\hspace*{5pt} We have just completed the proof of Theorem 7.
\section {$Spin(12)-$generators of the group $E_7$}
\hspace*{5pt} The simply connected compact Lie group $E_7$ is given by
$$
     E_7 = \{ \alpha \in \mbox{Iso}_C(\gP^C) \, | \, \alpha(P \times Q)\alpha^{-1} = \alpha P \times \alpha Q, <\alpha P, \alpha Q> = <P, Q> \}. $$
\noindent The group $E_7$ has subgroups
$$
    Spin_k(12) = \{ \alpha \in E_7 \, | \, \alpha\kappa_k = \kappa_k\alpha, \alpha\mu_k = \mu_k\alpha \} \quad(k = 1, 2) $$
\noindent where $\kappa_k$ and $\mu_k$ are defined by
\begin{center}
$\begin{array}{l}
   \kappa_k(X, Y, \xi, \eta) = (-(E_k, X)E_k + 4E_k \times (E_k \times X), (E_k, Y)E_k - 4E_k \times (E_k \times Y), - \xi, \eta), \\
   \mu_k(X, Y, \xi, \eta) = (2E_k \times Y + \eta E_k, 2E_k \times X + \xi E_k, (E_k, Y), (E_k, X)), 
\end{array} $
\end{center}
\noindent respectively, {\it e.g.}, when $k = 1$, for $P = ((\xi_1, \xi_2, \xi_3; x_1, x_2, x_3), (\eta_1, \eta_2, \eta_3; y_1, y_2, y_3), \xi, \eta)\in \gP^C$,\begin{center}
$\begin{array}{l}
   \kappa_1P = ((-\xi_1, \xi_2, \xi_3; x_1, 0, 0), (\eta_1, -\eta_2, -\eta_3; -y_1, 0, 0), -\xi, \eta), \\
   \mu_1P = ((\eta, \eta_3, \eta_2; - y_1, 0, 0), (\xi, \xi_3, \xi_2; -x_1, 0, 0), \eta_1, \xi_1). 
\end{array} $
\end{center} 
\noindent Then $Spin_k(12)$ is isomorphic to the usual spinor group $Spin(12)$ ([2], [4]).
\vspace{2mm}

\hspace*{5pt} {\bf Lemma 8.} {\it For an element}
      $P = ((\xi_1, \xi_2, \xi_3; x_1, x_2, x_3), 
           (\eta_1, \eta_2, \eta_3; y_1, y_2, y_3), \xi, \eta) 
\in \gP^C$ {\it satisfying} $P \times P = 0$, {\it it holds the following}
$$
   (1) \quad \xi_1\eta_1 + \xi_2\eta_2 + \xi_3\eta_3 + 2(x_1, y_1) + 2(x_2, y_2) + 2(x_3, y_3) - 3\xi\eta = 0, $$
$$  
    \begin{array}{clcl}
    (2) & \xi_2\xi_3 - \eta_1\eta - x_1\overline{x_1} = 0, &
    \quad (3) & \xi_3\xi_1 - \eta_2\eta - x_2\overline{x_2} = 0, \\
    (4) & \xi_1\xi_2 - \eta_3\eta - x_3\overline{x_3} = 0, &
    \quad(5) & \xi_1x_1 + \eta y_1 - \overline{x_2x_3} = 0, \\
    (6) & \xi_2x_2 + \eta y_2 - \overline{x_3x_1} = 0, &
    \quad(7) & \xi_3x_3 + \eta y_3 - \overline{x_1x_2} = 0, \\
    (8) & \eta_2\eta_3 - \xi_1\xi - y_1\overline{y_1} = 0, &
    \quad(9) & \eta_3\eta_1 - \xi_2\xi - y_2\overline{y_2} = 0, \\ 
    (10) & \eta_1\eta_2 - \xi_3\xi - y_3\overline{y_3} = 0, &
    \quad(11) & \eta_1y_1 + \xi x_1 - \overline{y_2y_3} = 0, \\
    (12) & \eta_2y_2 + \xi x_2 - \overline{y_3y_1} = 0, &
    \quad(13) & \eta_3y_3 + \xi x_3 - \overline{y_1y_2} = 0, \\
    (14) & \eta_3x_1 + \xi_2y_1 + \overline{y_2x_3} = 0, &
    \quad(15) & \eta_3x_2 + \xi_1y_2 + \overline{x_3y_1} = 0, \\
    (16) & \eta_2x_3 + \xi_1y_3 + \overline{y_1x_2} = 0, &
    \quad(17) & \eta_1x_3 + \xi_2y_3 + \overline{x_1y_2} = 0.   
\end{array}   $$
\vspace{2mm}

\hspace*{5pt} {\bf Proof} These are immediate from the straightforward computation of $P \times P = 0$. (Note that those are not all of the relations followed by $P \times P = 0$.)
\vspace{3mm}

\hspace*{5pt} {\bf Lemma 9.} (1) {\it For any element $P \in \gP^C$, there exists some element $\alpha_1 \in Spin_1(12)$ such that}
$$
   \alpha_1P = ((\xi_1, 0, 0; 0, x_2, x_3), 
           (\eta_1, \eta_2, \eta_3; 0, y_2, y_3), \xi, \eta). $$
\noindent {\it In particular, if an element} $P = ((0,\xi_2, \xi_3; x_1, 0, 0), 
           (\eta_1, 0, 0; 0, 0, 0), 0, \eta) \in \gP^C$ {\it satisfies the conditions $P \times P = 0$ and $<P, P> = 1$, then there exists some element $\alpha_1 \in Spin_1(12)$ such that} 
$$
       \alpha_1P = \d{1}, \quad where \, \, \d{1} = (0, 0, 0, 1) \in \gP^C. $$

\hspace*{5pt} (2) {\it For any element $P \in \gP^C$, there exists some element $\alpha_2 \in Spin_2(12)$ such that 
$$
    \alpha_2P = ((0, \xi_2, 0; x_1, 0, x_3), 
           (\eta_1, \eta_2, \eta_3; y_1, 0, y_3), \xi, \eta). $$
\noindent In particular, if an element $P = ((\xi_1, 0, \xi_3; 0, x_2, 0), (0, \eta_2, 0; 0, 0, 0), 0, \eta) \in \gP^C$ satisfies the conditions $P \times P = 0$ and $<P, P> = 1$, then there exists some element $\alpha_2 \in Spin_2(12)$ such that}
$$
               \alpha_2P = \d{1}. $$

\hspace*{5pt} {\bf Proof} (1) The first half is the very [2] Proposition 4.(2). We shall now prove the latter half. For an element $P = ((0, \xi_2,        \xi_3; x_1, 0, 0), (\eta_1, 0, 0; 0, 0, 0), 0, \eta) \in \gP^C$, act $\alpha_1 \in Spin_1(12)$ that is given in the first half which is composed of the elements of $Spin_1(12)$ defined in Lemmas 2, 3 and 4, on $P$. Then we get
$$
   \alpha_1P = ((0, 0, 0; 0, 0, 0), ({\eta_1}', 0, 0; 0, 0, 0), 0, {\eta}') 
                 = P_1, $$
\noindent noting that the subspaces $<\gP^C>_1, {<\gP^C>_1}'$ and ${<\gP^C>_1}''$ of $\gP^C$ are invariant under the action of the elements of $Spin_1(12)$ defined in Lemmas 2, 3 and 4, respectively, where
$$
      <\gP^C>_1 = \{((\xi_1, 0, 0; 0, 0, 0), 
           (0,\eta_2, \eta_3; y_1, 0, 0), \xi, 0) \in \gP^C \}, $$
$$
      {<\gP^C>_1}' = \{((0, \xi_2, \xi_3; x_1, 0, 0), 
           (\eta_1, 0, 0; 0, 0, 0), 0, \eta) \in \gP^C \}, $$
$$      {<\gP^C>_1}'' = \{((0, 0, 0; 0, x_2, x_3), 
           (0, 0, 0; 0, y_2, y_3), 0, 0) \in \gP^C \}. $$
\noindent From $P \times P = 0$, we have ${\eta_1}'{\eta}' = 0$ by Lemma 8.(2). As a result, the argument is devided into the following three cases:
\vspace{1mm}

\hspace*{5pt} (I) Case ${\eta_1}' = 0, \eta' \neq 0.  P_1$ is of the form $P_1 = ((0, 0, 0; 0, 0, 0), (0, 0, 0; 0, 0, 0), 0, {\eta}')$. Now, for $\theta \in C$ satifying $(\tau\theta)\theta = 1$, define the mapping $\epsilon_1(\theta) : \gP^C \to \gP^C$ as follows.
\begin{center}
$\begin{array}{l}
    \epsilon_1(\theta)((\xi_1, \xi_2, \xi_3; x_1, x_2, x_3), 
           (\eta_1, \eta_2, \eta_3; y_1, y_2, y_3), \xi, \eta)   \\
 \quad    = ((\theta^{-2}\xi_1, \xi_2, \xi_3; x_1, \theta^{-1}x_2, \theta^{-1}x_3), 
           (\theta^2\eta_1, \eta_2, \eta_3; y_1, \theta y_2, \theta y_3),
             \theta^2\xi, \theta^{-2}\eta). $$ 
\end{array} $
\end{center}
\noindent Then $\epsilon_1(\theta) \in Spin_1(12)$. Therefore, noting that $(\tau\eta')\eta' = <P_1, P_1> = 1,$ choose $\theta \in C$ such that $\theta^2 = \eta'$ and set $\epsilon_1(\theta)$. Then we get $\epsilon_1(\theta)P_1 = \d{1}$.
\vspace{1mm}

\hspace{5pt} (II) Case ${\eta_1}' \neq 0, {\eta}' = 0.$ By considering $\gamma_1(\pi/2)P_1$, where $\gamma_1(\pi/2) \in Spin_1(12)$ of Lemma 4.(1), this can be reduced to Case (I).
\vspace{1mm} 
 
\hspace*{5pt} (III) Case ${\eta_1}' = \eta' = 0.$ This does not occur, because $<P_1, P_1> = 1.$
\vspace{1mm}

\hspace{5pt} (2)  It is similarly verified by using $Spin_2(12)$ instead of $Spin_1(12)$ in the proof of (1).
\vspace{3mm}

\hspace*{5pt} {\bf Theorem 10.} {\it Any element $\alpha \in E_7$ can be represented as}
$$
    \alpha = \alpha_1\alpha_2{\alpha_1}'{\alpha_2}'{\alpha_1}'', \quad
       \alpha_1, {\alpha_1}', {\alpha_1}'' \in Spin_1(12),
       \alpha_2, {\alpha_2}'  \in Spin_2(12). $$

\hspace*{5pt} {\bf Proof} For a given element $\alpha \in E_7$, it suffices to show that there exist $\alpha_1, {\alpha_1}' \in Spin_1(12)$ and $\alpha_2 \in Spin_2(12)$ such that ${\alpha_1}'\alpha_2\alpha_1\alpha\d{1} = \d{1}$. In fact, since an element $\alpha \in E_7$ belongs to $E_6 (\subset E_7)$ if and only if $\alpha$ fixes an alement $\d{1}$, {\it i.e.,} $\alpha\d{1} = \d{1}$ ([4]), it follows ${\alpha_1}'\alpha_2\alpha_1\alpha \in E_6$, which implies that $\alpha \in E_7$ can be represented as a required form by Theorem 7. Now, set
$$
    \alpha\d{1} = ((\xi_1, \xi_2, \xi_3; x_1, x_2, x_3), 
           (\eta_1, \eta_2, \eta_3; y_1, y_2, y_3), \xi, \eta)    
                 = P_0 \in \gP^C. $$
\noindent Then, by Lemma 9.(1), we can take $\alpha_1 \in Spin_1(12)$ such that
 $$
    \alpha_1P_0 = (({\xi_1}', 0, 0; 0, {x_2}', {x_3}'), 
    ({\eta_1}', {\eta_2}', {\eta_3}'; 0, {y_2}', {y_3}'), {\xi}', {\eta}')    
                 = P_1. $$
\noindent From $P_1 \times P_1 = 0,$ we have ${\eta_1}'\eta' = 0$ by Lemma 8.(2). As a result, the argument is devided into the following three cases :
\vspace{1mm}

\hspace*{5pt} (I) Case    ${\eta_1}' = 0, \eta' \neq 0$. By Lemma 8.(6) and (7), we get ${y_2}' = {y_3}' = 0$. Furthermore we get $\xi' = 0$ by Lemma 8.(1). Therefore $P_1$ is of the form
$$
     P_1 = (({\xi_1}', 0, 0; 0, {x_2}', {x_3}'), 
    (0, {\eta_2}', {\eta_3}'; 0, 0, 0), 0, {\eta}'). $$  
\noindent  Then, by Lemma 8.(8), we have ${\eta_2}'{\eta_3}' = 0$. Hence there are three cases to be considered.
\vspace{1mm}

\hspace*{5pt} (I.A)  Case    ${\eta_2}' = 0, {\eta_3}' \neq 0$. By Lemma 8.(15), we get ${x_2}' = 0$, that is, $P_1$ is of the form
$$
     P_1 = (({\xi_1}', 0, 0; 0, 0, {x_3}'), 
    (0, 0, {\eta_3}'; 0, 0, 0), 0, {\eta}'). $$ 
\noindent Then, applying Lemma 9.(2) to $\alpha_1(\pi/2)P_1$, where $\alpha_1(\pi/2) \in Spin_1(9),$ we can obtain that there exists some element   $\alpha_2 \in Spin_2(12)$ such that $\alpha_2\alpha_1(\pi/2)P_1 = \d{1}.$
 \vspace{1mm}

\hspace*{5pt} (I.B)  Case ${\eta_2}' \neq 0, {\eta_3}' = 0$. By Lemma 8.(16), we get ${x_3}' = 0$, that is, $P_1$ is of the form
$$
     P_1 = (({\xi_1}', 0, 0; 0, {x_2}', 0), 
    (0, {\eta_2}', 0; 0, 0, 0), 0, {\eta}'). $$ 
\noindent Thus we can easily obtain the required result by Lemma 9.(2).
\vspace{1mm}

\hspace*{5pt} (I.C)  Case ${\eta_2}' = {\eta_3}' = 0$. $P_1$ is of the form
$$
     P_1 = (({\xi_1}', 0, 0; 0, {x_2}' , {x_3}'), 
    (0, 0, 0; 0, 0, 0), 0, {\eta}'). $$ 
\noindent Here we distinguish the following cases :
\vspace{1mm}

\hspace*{5pt} (I.C.1) When ${x_2}' \neq 0, {x_3}' \neq 0$. By Lemma 6.(2), we can take ${\alpha_1}' \in Spin_1(9)$ such that
$$
     {\alpha_1}'P_1 = (({\xi_1}'', 0, 0; 0, {x_2}'', {x_3}''), 
    (0, 0, 0; 0, 0, 0), 0, \eta''), \quad  {\xi_1}'' = {\xi_1}', \eta'' = \eta' \in C, {x_2}'' \in {\CC}^C, {x_3}'' \in \CC. $$  
\noindent Then, by Lemma 8.(4) we have ${x_3}''\overline{{x_3}''} = 0$, hence ${x_3}'' = 0$. Thus we easily obtain the required result by Lemma 9.(2).
\vspace{1mm}

\hspace*{5pt} (I.C.2) When ${x_2}' = 0, {x_3}' \neq 0$. Considering $\alpha_1(\pi/2)P_1$, where $\alpha_1(\pi/2) \in Spin_1(9)$, we can easily obtain the required result by Lemma 9.(2). 
\vspace{1mm}

\hspace*{5pt} (I.C.3) When ${x_2}' \neq 0, {x_3}' = 0$. We can easily obtain the required result by Lemma 9.(2). 
\vspace{1mm}

\hspace*{5pt} (I.C.4) When ${x_2}' = {x_3}' = 0$. We can easily obtain the required result by Lemma 9.(2). 
\vspace{1mm}

\hspace*{5pt} (II) Case ${\eta_1}' \neq 0, \eta' = 0$. By considering $\delta_1(\pi/2)P_1$, where $\delta_1(\pi/2) \in Spin_1(12)$ of Lemma 4.(3), this can be reduced to Case (I).
\vspace{1mm}

\hspace*{5pt} (III) Case ${\eta_1}' = \eta' = 0.$  $P_1$ is of the form
$$            P_1 = (({\xi_1}', 0, 0; 0, {x_2}', {x_3}'), 
        (0, {\eta_2}', {\eta_3}'; 0, {y_2}', {y_3}'), \xi', 0). $$ 
\noindent Now, as is similar to Lemma 9.(1), we obtain that, for any element $P \in \gP^C$, there exists some element $\alpha_1 \in Spin_1(12)$ such that
$$
    \alpha_1P = ((\xi_1, \xi_2, \xi_3; 0, x_2, x_3), 
           (\eta_1, 0, 0; 0, y_2, y_3), \xi, \eta).  $$
\noindent Note that the invariant subspaces $<\gP^C>_1, {<\gP^C>_1}'$ and ${<\gP^C>_1}'' $ of $\gP^C$ under the action of the elements of $Spin_1(12)$ defined in Lemmas 2, 3 and 4. Then, applying the result above to the present case (III), we can take ${\alpha_1}' \in Spin_1(12)$ such that
$$
       {\alpha_1}'P_1 = (({\xi_1}'',0, 0; 0, {x_2}'', {x_3}''), 
           (0, 0, 0; 0, {y_2}'', {y_3}''), \xi'', 0) = P_2. $$
\noindent Therefore we have ${\xi_1}''\xi'' = 0$ by Lemma 8.(8). Hence there  are three cases to be considered.

\hspace*{5pt} (III.A) Case ${\xi_1}'' = 0, \xi'' \neq 0.$  By Lemma 8.(12) and (13), we get ${x_2}'' = {x_3}'' = 0$. Then $P_2$ is of the form
$$
       P_2 = ((0, 0, 0; 0, 0, 0), 
           (0, 0, 0; 0, {y_2}'', {y_3}''), \xi'', 0). $$
\noindent Thus, by considering $\gamma_1(\pi/2)P_2$, where $\gamma_1(\pi/2) \in Spin_1(12)$ of Lemma 4.(1), this can be reduced to Case (I.C).
\vspace{1mm}

\hspace*{5pt} (III.B) Case ${\xi_1}'' \neq 0, \xi'' = 0$. By Lemma 8.(15) and (16), we get ${y_2}'' = {y_3}'' = 0.$ Therefore this is reduced to Case (I.C).
\vspace{1mm}

\hspace*{5pt} (III.C) Case ${\xi_1}'' = \xi'' = 0$. $P_2$ is of the form
$$
     P_2 = ((0,0, 0; 0, {x_2}'', {x_3}''), 
           (0, 0, 0; 0, {y_2}'', {y_3}''), 0, 0). $$
\noindent Here we distinguish the following cases :
\vspace{1mm}

\hspace*{5pt} (III.C.1) When ${x_2}'' \neq 0.$  By Lemma 9.(2), there exists some element $\alpha_2 \in Spin_2(12)$ such that
$$
   \alpha_2P_2 = ((0, \xi_2^{(3)}, 0; x_1^{(3)}, 0, x_3^{(3)}),
        (\eta_1^{(3)}, \eta_2^{(3)}, \eta_3^{(3)}; y_1^{(3)}, 0, y_3^{(3)}),
                   \xi^{(3)}, \eta^{(3)}) = P_3. $$
\noindent Here, by Lemma 8.(3), we have $\eta_2^{(3)}\eta^{(3)} = 0$. Hence there are three cases to be considered.
\vspace{1mm}

\hspace*{5pt} (III.C.1.1) Case $\eta_2^{(3)} = 0, \eta^{(3)} \neq 0.$ By Lemma 8.(5) and (7), we get $y_1^{(3)} = y_3^{(3)} = 0$.  Furthermore, we get $\xi^{(3)} = 0$ by Lemma 8.(1). Then $P_3$ is of the form
$$
      P_3 = ((0, \xi_2^{(3)}, 0; x_1^{(3)}, 0, x_3^{(3)}),
                  (\eta_1^{(3)}, 0, \eta_3^{(3)}; 0, 0, 0),
                   0, \eta^{(3)}). $$
\noindent Here, by lemma 8.(9), we have $\eta_3^{(3)}\eta_1^{(3)} = 0.$ Hence there are three cases to be considered.
\vspace{1mm}

\hspace*{5pt} (III.C.1.1.1) Case $\eta_1^{(3)} = 0, \eta_3^{(3)} \neq 0.$  By Lemma 8.(14), we get $x_1^{(3)} = 0$. Then, considering $\alpha_2(\pi/2)P_3$, where $\alpha_2(\pi/2) \in Spin_2(9)$, we can easily obtain the required result by Lemma 9.(1).
\vspace{1mm}

\hspace*{5pt} (III.C.1.1.2) Case $\eta_1^{(3)} \neq 0, \eta_3^{(3)} = 0$. By Lemma 8.(17), we get $x_3^{(3)} = 0$. Then we can easily obtain the required result by Lemma 9.(1).
\vspace{1mm}

\hspace*{5pt}  (III.C.1.1.3) Case $\eta_1^{(3)} = \eta_3^{(3)} = 0$. $P_3$ is of the form
$$
      P_3 = ((0, \xi_2^{(3)}, 0; x_1^{(3)}, 0, x_3^{(3)}),
                  (0, 0, 0; 0, 0, 0),
                   0, \eta^{(3)}). $$
\noindent Here we distinguish the following cases :
\vspace{1mm}

\hspace*{5pt} (III.C.1.1.3.(i)) When $x_1^{(3)} \neq 0, x_3^{(3)} \neq 0$. As is similar to Lemma 6.(2), we obtain that there exists some element ${\alpha_2}' \in Spin_2(9)$ such that
\begin{center}
$
      {\alpha_2}'P_3 = ((0, \xi_2^{(4)}, 0; x_1^{(4)}, 0, x_3^{(4)}),
                  (0, 0, 0; 0, 0, 0),
                   0, \eta^{(4)}) = P_4$, \quad 
  $\begin{array}{l}
                \xi_2^{(4)} = \xi_2^{(3)},  \eta^{(4)} = \eta^{(3)} \in C, \\
                 x_1^{(4)} \in \CC , x_3^{(4)} \in \CC^C. 
   \end{array} $
   \end{center}
   \noindent Then, by Lemma 8.(2), we have $x_1^{(4)}\overline{x_1^{(4)}} = 0$,  hence  $x_1^{(4)} = 0$. Thus, considering $\alpha_2(\pi/2)P_4$, where $\alpha_2(\pi/2) \in Spin_2(9)$, we can easily obtain the required result by Lemma 9.(1).
\vspace{1mm}

\hspace*{5pt} (III.C.1.1.3.(ii)) When $x_1^{(3)} = 0, x_3^{(3)} \neq 0$.  Considering $\alpha_2(\pi/2)P_3$, where $\alpha_2(\pi/2) \in Spin_2(9)$, we can easily obtain the required result by Lemma 9.(1).
\vspace{1mm}

\hspace*{5pt} (III.C.1.1.3.(iii)) When $x_1^{(3)} \neq 0, x_3^{(3)} = 0$. We easily obtain the required result by Lemma 9.(1).
\vspace{1mm}

\hspace*{5pt} (III.C.1.1.3.(iv)) When $x_1^{(3)} = x_3^{(3)} = 0$. We easily obtain the required result by Lemma 9.(1).
\vspace{1mm}

\hspace*{5pt} (III.C.1.2) Case $\eta_2^{(3)} \neq 0, \eta^{(3)} = 0$. By considering $\gamma_2(\pi/2)P_3$, where $\gamma_2(\pi/2) \in Spin_2(12)$ of Lemma 4.(2), this can be reduced to Case (III.C.1.1). 
\vspace{1mm}

\hspace*{5pt} (III.C.1.3) Case $\eta_2^{(3)} = \eta^{(3)} = 0$. This does not occur. In fact, note that the subspace $<\gP^C>_2$ of $\gP^C$ is invariant under the action of the elements of $Spin_2(12)$ defined in Lemmas 2, 3 and 4, where $<\gP^C>_2 = \{((\xi_1, 0, \xi_2,; 0, x_2, 0), (0, \eta_2, 0; 0, 0, 0), 0, \eta) \in \gP^C \}$. Then, for $P_3 = \alpha_2P_2$, that is,
$$
     ((0, \xi_2^{(3)}, 0; x_1^{(3)}, 0, x_3^{(3)}),
       (\eta_1^{(3)}, 0, \eta_3^{(3)}; y_1^{(3)}, 0, y_3^{(3)}),
                   \xi^{(3)}, 0) $$
 $$         = \alpha_2((0, 0, 0; 0, {x_2}'', {x_3}''),
                   (0, 0, 0; 0, {y_2}'', {y_3}''), 0, 0), $$
\noindent where $\alpha_2 \in Spin_2(12)$, the condition $\eta_2^{(3)} = \eta^{(3)} = 0$ contradicts $ {x_2}'' \neq 0$.
\vspace{1mm}

\hspace*{5pt} (III.C.2) When ${x_2}'' = 0, {x_3}'' \neq 0.$  By considering $\alpha_1(\pi/2)P_2$, where $\alpha_1(\pi/2) \in Spin_1(9)$, this can be reduced to Case (III.C.1). 
\vspace{1mm}

\hspace*{5pt} (III.C.3) When ${x_2}'' = {x_3}'' = 0, {y_3}'' \neq 0.$  By considering $\gamma_1(\pi/2)P_2$, where $\gamma_1(\pi/2) \in Spin_1(12)$, this can be reduced to Case (III.C.1). 
\vspace{1mm}

\hspace*{5pt} (III.C.4) When ${x_2}'' = {x_3}'' = {y_3}'' = 0, {y_2}'' \neq 0.$  By considering $\alpha_1(\pi/2)P_2$, where $\alpha_1(\pi/2) \in Spin_1(9)$, this can be reduced to Case (III.C.3). 

\vspace{1mm}

\hspace*{5pt} (III.C.5) When ${x_2}'' = {x_3}'' = {y_2}'' = {y_3}'' = 0.$ It is obvious that this does not occur. 
\vspace{1mm}

\hspace*{5pt} We have just completed the proof of Theorem 10.     
\vspace{2mm} 

   \hspace*{5pt} {\bf Conjecture}  We know that the simply connected compact excetional Lie group $E_8$ has subgroups 
   $Ss_k(16) = (E_8)^{\sigma_k} $ (where $\sigma_k = \mbox{exp}\pi\kappa_k), k = 1, 2, 3$ 
   (which is isomorphic to $Spin(16)/\Z_2$ not $SO(16))$. Now the authors do not know if 
  $Ss_1(16)$ and $Ss_2(16)$ generate the group $E_8$?
    \vspace{4mm}

\centerline {\bf References}
\vspace{2mm}

[1] C. Chevalley, {\it Theorey of Lie groups}, Princeton Univ. Press, New Jersey, 1946.
\vspace{1mm}

[2] T. Miyasaka, O. Yasukura and I. Yokota, {\it Diagonalization of an element $P$ of $\gP^C$ by the compact Lie group $E_7$}, Tsukuba J.Math., 22(1998), 687-703.
\vspace{1mm}

[3] I. Yokota, {\it Realization of involutive automorphisms $\sigma$ and $G^{\sigma}$ of exceptional linear Lie groups $G,$} Part I, $G = G_2, F_4$ {\it and} $E_6$, Tsukuba J. Math., 14(1990), 185-223.
\vspace{1mm}

[4] I. Yokota, {\it Realization of involutive automorphisms $\sigma$ and $G^{\sigma}$ of exceptional linear Lie groups $G,$} Part II, $G = E_7$, Tsukuba J. Math., 14(1990), 379-404.

\end{document}